\documentclass[12pt]{amsart}
\usepackage{amsthm, amsmath, amssymb, enumerate, bm}
\newtheorem{thm}{Theorem}[section]
\newtheorem{prop}[thm]{Proposition}
\newtheorem{lemm}[thm]{Lemma}

\theoremstyle{definition}

\newcommand{\C}{\mathbb{C}}
\newcommand{\R}{\mathbb{R}}
\newcommand{\Q}{\mathbb{Q}}
\newcommand{\Z}{\mathbb{Z}}
\newcommand{\BH}{\mathbb{H}}
\newcommand{\CA}{\mathcal{A}}
\newcommand{\CB}{\mathcal{B}}

\newcommand{\CI}{\mathcal{I}}
\newcommand{\CJ}{\mathcal{J}}
\newcommand{\CJF}{\mathcal{J} (\F)}
\newcommand{\RI}{\mathrm{i}}
\newcommand{\RSYM}{\mathrm{sym}}
\newcommand{\F}{F}
\newcommand{\FG}{F[G]}
\newcommand{\JS}{(X,S)}
\newcommand{\ZA}{Z(\mathcal{A})}
\newcommand{\CHAR}{\operatorname{char}}
\newcommand{\IRR}{\operatorname{Irr}}
\newcommand{\SPAN}{\operatorname{span}}
\newcommand{\SPIN}{\operatorname{JSpin}}
\newcommand{\TRACE}{\operatorname{tr}}

\begin{document}
\title[Adjacency algebras of thin Jordan schemes]
{Structure of adjacency algebras of thin Jordan schemes}

\author[Y.Ohisa]{Yudai Ohisa}
\address{Faculty of Science, Shinshu University, 
 3-1-1 Asahi, Matsumoto, Nagano 390-8621, Japan}
\email{yudai.ohisa@gmail.com}

\author[K.Yanagisawa]{Kyohei Yanagisawa}
\address{Faculty of Science, Shinshu University, 
 3-1-1 Asahi, Matsumoto, Nagano 390-8621, Japan}
\email{kyohei.yanagisawa@gmail.com}

\begin{abstract}
We prove that adjacency algebras of Jordan schemes over fields of characteristic zero are semisimple.
We then classify the adjacency algebras of thin Jordan schemes over the fields of real numbers or complex numbers.
The simple summands of the adjacency algebras over these fields are restricted to specific types of simple Jordan algebras.
\end{abstract}


\maketitle

\section{Introduction}
Association schemes describe partitions of the set of ordered pairs of a finite set into relations whose adjacency matrices satisfy prescribed product rules.
The entries of a product of two adjacency matrices count intermediate points connecting two points through the corresponding relations in a specified order.
The linear span of these matrices is closed under ordinary matrix product and forms an associative algebra, called the adjacency algebra of the association scheme.
Studying its representations and its decomposition into simple algebras provides algebraic methods for investigating the structure of the association scheme.

Jordan schemes generalize association schemes by replacing the ordinary matrix product with the Jordan product $A \star B := \frac{AB + BA}{2}$ where $A, B$ are square matrices.
The linear span of the adjacency matrices of a Jordan scheme forms a Jordan algebra under the Jordan product.
This algebra is called the \emph{adjacency algebra} of a Jordan scheme.
The Jordan product is commutative but is generally non-associative.
It also preserves the space of symmetric matrices, which need not be closed under the ordinary matrix product.
For an association scheme, merging each relation with its transpose gives adjacency matrices that span the symmetric part of its adjacency algebra.
This subspace is closed under the Jordan product, so the resulting structure is a Jordan scheme.
However, it need not be closed under the ordinary matrix product, and hence the resulting Jordan scheme need not be an association scheme.

Adjacency algebras of Jordan schemes may similarly provide natural algebraic methods for investigating the structure of Jordan schemes.
Since these algebras are generally non-associative, we study them using the structure theory of Jordan algebras.
Two natural questions are whether they are semisimple and which simple Jordan algebras occur in their decompositions.

Our first result applies to arbitrary Jordan schemes.
It is known that the adjacency algebra of a symmetric Jordan scheme over $\R$ is semisimple \cite{MR4462635}.
We extend this result to arbitrary Jordan schemes over fields of characteristic $0$.
Consequently, these algebras decompose into direct sums of simple Jordan algebras.

Hanaki and Yoshikawa constructed proper Jordan schemes related to quaternion and octonion algebras and determined the structures of their adjacency algebras over the field of real numbers \cite{arxiv2509.01865}.
Muzychuk, Pech, and Woldar classified thin Jordan schemes \cite{MR5094503}.
Their classification separates thin Jordan schemes into non-regular and regular cases.
The regular case further divides into thin association schemes corresponding to finite groups and regular thin Jordan schemes corresponding to non-associative RA-loops.

We then determine these decompositions over $\R$ or $\C$ for thin Jordan schemes.
In the non-regular case, the simple summands over $\C$ are spin factors with two-dimensional vector space parts, and over $\R$, they are spin factors of signature $(1, 1)$ or realifications of the complex spin factors.
In the regular case corresponding to finite groups, namely the thin association schemes case, the decomposition is obtained from the classical Wedderburn decompositions of group algebras.
In the regular case corresponding to non-associative RA-loops, the simple summands over $\C$ are copies of $\C$ or spin factors with seven-dimensional vector space parts, and over $\R$, they are copies of $\R$, copies of $\C$ regarded as real Jordan algebras, spin factors of signatures $(0, 7)$ and $(4, 3)$, or realifications of the complex spin factors.

This paper is organized as follows.
Section 2 recalls the definitions and the classification of thin Jordan schemes, together with the required facts about Jordan algebras.
Section 3 proves semisimplicity and determines the simple Jordan algebra decompositions in the non-regular case and the regular cases corresponding to finite groups and non-associative RA-loops.

\section{Preliminaries}
Unless otherwise stated, let $\F$ be $\R$ or $\C$.
We denote the identity matrix by $I$, the all-ones matrix by $J$, and the zero matrix by $O$.
The transpose of a matrix $M$ is denoted by ${}^t\! M$.
For matrices $A, B$, we denote the ordinary matrix product by $AB$.

\subsection{Jordan scheme}
Let $X$ be a finite set, and let $M_X(\F)$ denote the full matrix algebra over $\F$ whose rows and columns are indexed by $X$.
We consider a set of $\{0, 1\}$-matrices $S = \{I = A_0, A_1, \ldots, A_d\}\subset M_X(\F)$ satisfying the following conditions:
\begin{enumerate}
\item[(JS1)] $\sum_{i = 0}^{d}A_i = J$, and $S$ contains no zero matrix, 
\item[(JS2)] For $i\in\{0, 1, \ldots, d\}$, there exists $j\in\{0, 1, \ldots, d\}$ such that ${}^tA_i = A_j$, 
\item[(JS3)] For $i, j, k\in\{0, 1, \ldots, d\}$, there exists a rational number $p_{ij}^{k}\in\frac{1}{2}\Z$ such that $A_i\star A_j = \sum_{k = 0}^{d}p_{ij}^{k}A_k$.
\end{enumerate}
Then the pair $(X, S)$ is called a \emph{Jordan scheme}.

A Jordan scheme is called \emph{thin} if, for every $A_i\in S$, each row of $A_i$ contains at most one entry equal to $1$.
A Jordan scheme is called \emph{regular} if, for every $A_i\in S$, the number of entries equal to $1$ in each row of $A_i$ is independent of the choice of the row.

\subsection{Jordan algebra}\label{subsec:jordan-algebra}
Let $\F$ be a field of $\CHAR(\F)\neq 2$.
A vector space $V$ over $\F$ equipped with a bilinear product $\star$ is called a \emph{Jordan algebra} over $\F$ if all $x, y\in V$ satisfy
\begin{enumerate}
\item[(JA1)] commutativity:  $x \star y = y \star x$, 
\item[(JA2)] Jordan identity: $((x \star x) \star y) \star x = (x \star x) \star (y \star x)$.
\end{enumerate}

By condition (JS3), the $\F$-vector space $\F S:=\bigoplus_{i=0}^{d}\F A_i$ is closed under the Jordan product.
Since the Jordan product is commutative and satisfies the Jordan identity, $\F S$ is a Jordan algebra.
This Jordan algebra is called the \emph{adjacency algebra} of a Jordan scheme $(X, S)$ over $\F$.

Let $\CA$ be a finite-dimensional Jordan algebra over $\F$.
A subset $\CI$ of $\CA$ is called an \emph{ideal} of $\CA$ if
\begin{enumerate}
\item[(I1)] $\CI$ is a subspace of the vector space $\CA$, 
\item[(I2)] $a \star x \in \CI$ for all $a \in \CA, x \in \CI$.
\end{enumerate}
We write $a^{\star 2} = a\star a$ for the Jordan product square of $a$.
For matrices, the Jordan powers agree with the ordinary matrix powers.
An element $a \in \CA$ is called \emph{nilpotent} if $a^{\star r}=0$ for some positive integer $r$, where $a^{\star r}$ denotes the $r$-th power of $a$ with respect to the Jordan product.
An ideal of $\CA$ is called \emph{nil} if all of its elements are nilpotent.

The Jordan algebra $\CA$ is called \emph{simple} if it contains a non-nilpotent element and has no ideals other than $\{0\}$ and $\CA$.
For the classification of finite-dimensional simple Jordan algebras over $\C$, see \cite[Renaissance Structure Theorem, I.~2.13, p.~68]{MR2014924}.
For the classification over $\R$, see \cite[Lemma~3.7 and Theorems~5.9 and~5.10]{MR3391193}.

The Jordan algebra $\CA$ is called \emph{semisimple} if it contains no nonzero nil ideal.
By \cite[Theorem 11, p.~65]{MR1718170}, if $\CA$ is semisimple then it is a direct sum of simple Jordan algebras and these summands are exactly the minimal nonzero ideals of $\CA$.
Namely
\[
\CA = \bigoplus_{i = 1}^{k}\CA_i, \quad \text{where each $\CA_i$ is a simple Jordan algebra of $\CA$.}
\]
Note that this decomposition is unique up to a permutation of the summands.

For each $a\in\CA$, let $T_a\colon\CA\to\CA$ be the linear operator defined by $T_a(x) = a\star x$ for $x\in\CA$.
For a finite-dimensional Jordan algebra $\CA$ over a field $\F$ of characteristic $0$, define its trace form by
\[
\tau_{\F}(a, b):=\TRACE\left(T_{a\star b}\right)\quad(a, b\in\CA), 
\]
where the trace is taken on $\CA$.
Then $\CA$ is semisimple if and only if $\tau_{\F}$ is nondegenerate \cite[Chapter~III, \S~2, Theorem~6, p.~60]{MR1718170}.

The center of $\CA$ is defined by $\ZA:= \{a\in\CA\mid [T_a, T_b] = 0\text{ for all }b\in\CA\}$, where $[T_a, T_b]:= T_aT_b-T_bT_a$.
Since the Jordan product is commutative, for $a\in\ZA$ and $b, c\in\CA$ we have
\[
(T_{a\star b}-T_aT_b)c=(a\star b)\star c-a\star(b\star c)=-[T_a, T_c]b=0.
\]
Thus, $T_{a\star b}=T_aT_b$ for all $a\in\ZA$ and $b\in\CA$.

If $\CA$ is a Jordan subalgebra of a matrix algebra, then any $M\in\CA$ that commutes with every element of $\CA$ under the ordinary matrix product belongs to $\ZA$.
Indeed, for $A, B\in\CA$, we have
\[
(M\star A)\star B-M\star(A\star B) = \frac{1}{4}\left(AMB+BMA-MBA-ABM\right) = O.
\]

If $\CA$ is semisimple, then it has an identity element $1_{\CA}$, and its center $\ZA$, equipped with the restriction of the Jordan product $\star$, is a finite-dimensional commutative associative algebra with the same identity.
The decomposition $\CA = \bigoplus_{i = 1}^{k}\CA_i$ into simple ideals induces the decomposition $\ZA = \bigoplus_{i = 1}^{k}Z_i$, where $Z_i = \ZA\cap\CA_i$ is a field.
Let $e_i$ be the identity element of $Z_i$.
Then $e_1, \ldots, e_k$ are pairwise orthogonal primitive idempotents of $\ZA$ satisfying $1_{\CA} = e_1+\cdots+e_k$.
Here $\CA \star e_i:=\{a\star e_i\mid a\in\CA\}$.
Moreover, $\CA_i = \CA \star e_i$, and $e_i$ is the identity element of $\CA_i$.

Consequently, to determine the decomposition of a semisimple Jordan algebra $\CA$ into simple ideals, it suffices to decompose the identity into pairwise orthogonal primitive idempotents in its center $\ZA$ and form the corresponding ideals $\CA \star e_i$.

\begin{prop}\label{prop:realification}
Let $\CA$ be a unital simple Jordan algebra over $\C$.
Let $\CA_{\R}$ denote the realification of $\CA$, namely the Jordan algebra obtained by regarding $\CA$ as an $\R$-vector space with the same Jordan product.
Then $\CA_{\R}$ is simple.
\end{prop}
\begin{proof}
Let $\CI$ be a nonzero ideal of $\CA_{\R}$.
Since $\CA_{\R}$ and $\CA$ have the same underlying set, we have $z1_{\CA}\in\CA_{\R}$ for every $z\in\C$.
Hence, for $x\in\CI$, $zx = x \star (z1_{\CA}) \in \CI$.
Thus, $\CI$ is closed under product by complex scalars and is therefore a complex ideal of $\CA$.
Since $\CA$ is simple and $\CI\neq0$, we have $\CI=\CA_{\R}$.
\end{proof}

For a real Jordan algebra $\CA$, let $\CA_{\C} := \CA\otimes_{\R}\C$ denote its complexification.

\subsection{Spin factor}
Let $V$ be an $n$-dimensional vector space over $\F$ with a nondegenerate symmetric bilinear form $f$.
Define a bilinear product $\star$ on $\F \oplus V$ by 
\[(\alpha, x)\star(\beta, y)=(\alpha\beta+f(x, y), \alpha y+\beta x).\]
Since this product is commutative and satisfies the Jordan identity, $\F \oplus V$ is a Jordan algebra.
This Jordan algebra is called the \emph{Jordan algebra of a nondegenerate symmetric bilinear form} in \cite[I.~5]{MR251099}, and called the \emph{spin factor} in \cite[I.~1.9, p.~47]{MR2014924}.
In this paper, we use the term \emph{spin factor}.
If $\dim V > 1$, then it is simple \cite[Second Structure Theorem, IV.~5]{MR251099}.

By Sylvester's law of inertia, if $\F = \R$, then there exist nonnegative integers $p$ and $q$ with $p+q = n$ and a basis $\{v_1, \ldots, v_{p+q}\}$ of $V$ such that
\[
f(\bm{x}, \bm{y}) = \sum_{j = 1}^{p}x_jy_j - \sum_{j = p+1}^{p+q}x_jy_j \quad (\bm{x} = \sum_{j = 1}^{p+q}x_jv_j, \quad \bm{y} = \sum_{j = 1}^{p+q}y_jv_j).
\]
The pair $(p, q)$ is uniquely determined by $f$ and is called the \emph{signature} of $f$.
Thus, spin factors are classified by their signatures $(p, q)$.
Throughout, we denote the spin factor over $\R$ by $\SPIN_{\R} (p, q)$.
Similarly, over $\C$, every nondegenerate symmetric bilinear form is equivalent to the standard bilinear form.
Hence, complex spin factors are classified by the dimension $n$ of $V$.
Throughout, we denote the spin factor over $\C$ by $\SPIN_{\C} (n)$.
For $n>1$, Proposition~\ref{prop:realification} shows that $(\SPIN_{\C}(n))_{\R}$ is simple.

The product table for spin factors over $\R$ is given in TABLE~\ref{tab:real-spin-factor}, where $1:= (1, 0)$ and $\bm{v}_i:= (0, v_i)$.
For a complex spin factor, the basis of $V$ can be chosen so that $f(v_i, v_j) = \delta_{ij}$; hence, all diagonal entries $\bm{v}_i\star \bm{v}_i$ in the corresponding product table are equal to $1$.
The following product table will be used repeatedly in Section 3.
\begin{table}[htbp]
\caption{product table of a real spin factor of signature $(p, q)$}
\label{tab:real-spin-factor}
\centering
\renewcommand{\arraystretch}{1.5}
\begin{tabular}{c||c|ccc|ccc}
$\star$ & $1$ & $\bm{v}_1$ & $\cdots$ & $\bm{v}_p$ & $\bm{v}_{p+1}$ & $\cdots$ & $\bm{v}_{p+q}$ \\ \hline\hline
$1$ & $1$ & $\bm{v}_1$ & $\cdots$ & $\bm{v}_p$ & $\bm{v}_{p+1}$ & $\cdots$ & $\bm{v}_{p+q}$ \\ \hline
$\bm{v}_1$ & $\bm{v}_1$ & $+1$ & $\cdots$ & $0$ & $0$ & $\cdots$ & $0$ \\
$\vdots$ & $\vdots$ & $\vdots$ & $\ddots$ & $\vdots$ & $\vdots$ & $\ddots$ & $\vdots$ \\
$\bm{v}_p$ & $\bm{v}_p$ & $0$ & $\cdots$ & $+1$ & $0$ & $\cdots$ & $0$ \\ \hline
$\bm{v}_{p+1}$ & $\bm{v}_{p+1}$ & $0$ & $\cdots$ & $0$ & $-1$ & $\cdots$ & $0$ \\
$\vdots$ & $\vdots$ & $\vdots$ & $\ddots$ & $\vdots$ & $\vdots$ & $\ddots$ & $\vdots$ \\
$\bm{v}_{p+q}$ & $\bm{v}_{p+q}$ & $0$ & $\cdots$ & $0$ & $0$ & $\cdots$ & $-1$ \\
\end{tabular}
\end{table}

\subsection{Classification of non-regular thin Jordan schemes}\label{subsec:non-regular-thin-classification}
We recall the classification of non-regular thin Jordan schemes obtained by Muzychuk, Pech, and Woldar \cite{MR5094503}.

Let $G$ be a finite abelian group of order $n_G$ and let $\FG$ be the group algebra of $G$.
For each $g\in G$, let $P_g$ denote the permutation matrix of the left translation $h\mapsto gh$.
Let $\CA (G, G)$ denote the adjacency algebra of the thin association scheme associated with $G$.
The matrices $\{P_g\mid g\in G\}$ form a basis of $\CA (G, G)$ and satisfy $P_gP_h=P_{gh}$ for all $g, h\in G$.
Hence the assignment $g\mapsto P_g$ extends linearly to an $\F$-algebra isomorphism from $\FG$ onto $\CA (G, G)$.
Therefore we identify $\FG$ with $\CA (G, G)$ and denote it by $\FG$ throughout this paper.
We define
\[
\CJ(\FG) := \left\{\begin{bmatrix}A & B\\C & A\end{bmatrix}\;\middle|\; A, B, C\in\FG\right\} \subseteq M_{2n_G}(\F).
\]
Then $\CJ(\FG)$ is a Jordan algebra under the Jordan product.
Note that in \cite{MR5094503} the authors write $\CJ(G)$ for $\CJ(\CA(G, G))$.
Set
\[
I_2:=\begin{bmatrix}1&0\\0&1\end{bmatrix}, \quad E_{12}:=\begin{bmatrix}0&1\\0&0\end{bmatrix}, \quad E_{21}:=\begin{bmatrix}0&0\\1&0\end{bmatrix}, 
\]
\[
\CJF:=\left\{\begin{bmatrix}a&b\\c&a\end{bmatrix}\;\middle|\;a, b, c\in\F\right\}=\F I_2 \oplus \F E_{12} \oplus \F E_{21}\subseteq M_2(\F).
\]
Since $\CJF$ is closed under the Jordan product, $\CJF$ is a Jordan algebra.
Using the Kronecker product, every element of $\CJ(\FG)$ can be written uniquely as
\[
\begin{bmatrix}A&B\\C&A\end{bmatrix}=I_2\otimes A+E_{12}\otimes B+E_{21}\otimes C\quad(A, B, C\in\FG).
\]
Consequently, the Kronecker product induces an $\F$-linear isomorphism from $\CJF\otimes_{\F}\FG$ onto $\CJ(\FG)$.
Since $G$ is abelian, $\FG$ is commutative.
Thus, for $M, N\in\CJF$ and $A, B\in\FG$, we have
\[
\begin{aligned}
(M\otimes A)\star(N\otimes B)
&= \frac{(M\otimes A)(N\otimes B)+(N\otimes B)(M\otimes A)}{2}\\[4pt]
&= \frac{MN\otimes AB+NM\otimes BA}{2}\\[4pt]
&= (M\star N)\otimes AB.
\end{aligned}
\]
Hence the above linear isomorphism preserves the Jordan product, and therefore $\CJ(\FG)\cong \CJF \otimes_{\F}\FG$ as Jordan algebras.

\begin{prop}\label{prop:non-regular-classification}\cite[Proposition 3.5]{MR5094503}
Let $\JS$ be a non-regular thin Jordan scheme and let $\CA$ be its adjacency algebra.
Then there exists an abelian group $G$ and a permutation matrix $P$ such that $\CA = P^{-1}\CJ(\FG)P$.
\end{prop}
Conversely, for every finite abelian group $G$, $\CJ(\FG)$ is the adjacency algebra of a non-regular thin Jordan scheme \cite[Theorem 3.6]{MR5094503}.

\subsection{Classification of regular thin Jordan schemes}\label{subsec:regular-thin-classification}
We recall the classification of regular thin Jordan schemes obtained by Muzychuk, Pech, and Woldar \cite{MR5094503}.
We first recall a definition of an RA-loop (ring alternative loop).

A \emph{loop} is a nonempty set $L$ with a binary operation $\cdot$ satisfying the following conditions \cite{MR1433590}.
\begin{enumerate}
\item[(L1)] There exists $1_{L}\in L$ such that $1_{L}\cdot a=a\cdot 1_{L} = a$ for every $a\in L$.
\item[(L2)] For every $a, b\in L$, each of the equations $a\cdot x=b$ and $y\cdot a=b$ has a unique solution in $L$.
\end{enumerate}
For a commutative associative ring $R$, let $RL$ be the loop ring of $L$ over $R$.
The loop $L$ is called an \emph{RA-loop (ring alternative loop)} if $RL$ is alternative for every commutative associative ring $R$, namely, 
\[
(xx)y = x(xy) \quad \text{and}\quad (yx)x = y(xx) \quad \text{for all }x, y\in RL.
\]
Non-associative RA-loops are described by the following construction \cite{MR1351625}.
Let $G$ be a finite group satisfying $G/Z(G) \cong \Z_{2}^2$.
Then the commutator subgroup of $G$ is $G' = \{1_{G}, s\}\subseteq Z(G)$, where $1_{G}$ is the identity element of $G$.
Define an involution $g\mapsto g^*$ on $G$ by
\[
g^* = 
\begin{cases}
g, & \text{if }g\in Z(G), \\
gs, & \text{otherwise}.
\end{cases}
\]
Choose an arbitrary element $g_0\in Z(G)$ and define a binary operation $\diamond$ on the set $G\times\{0, 1\}$ as follows:
\[
(g, x) \diamond (h, y) = \left\{
\begin{array}{ll}
(gh, 0), & \text{if } x = 0, y = 0\\
(hg, 1), & \text{if } x = 0, y = 1\\
(gh^*, 1), & \text{if } x = 1, y = 0\\
(g_0h^*g, 0), & \text{if } x = 1, y = 1.
\end{array}
\right.
\]
Then $(G\times\{0, 1\}, \diamond)$ is an RA-loop, denoted by $L(G, *, g_0)$.

For a loop $(L, \diamond)$ and $a\in L$, let $\lambda_a$ denote the left translation $b\mapsto a\diamond b$ on $L$.
Let $\ell_a$ denote the permutation matrix of $\lambda_a$ with respect to a fixed ordering of $L$.
The columns correspond to the inputs of the left translation, and the rows correspond to their images.
For $v, w\in L$, the $(v, w)$-entry of $\ell_a$ is $1$ if $\lambda_a(w)=v$ and $0$ otherwise.
For $Y = \sum_{u\in L}c_u u\in\F L$, define $\ell_Y := \sum_{u\in L}c_u \ell_u$.
\begin{thm}\label{thm:regular-classification}
Let $(L, \diamond)$ be a finite loop.
The permutation matrices $\{\ell_a\}_{a\in L}$ of the left translations form a regular thin Jordan scheme if and only if $L$ is an RA-loop \cite[Theorem 3.17]{MR5094503}.
Moreover, every regular thin Jordan scheme is isomorphic to a scheme obtained in this way \cite[Proposition 3.9]{MR5094503}.
\end{thm}
If $L$ is a non-associative RA-loop, then $L$ is isomorphic to $L(G, *, g_0)$ for a suitable finite group $G$ and an element $g_0\in Z(G)$ as described above.
If $L$ is associative, then $L$ is a group and the resulting Jordan scheme is a thin association scheme.

\subsection{Group algebras}\label{subsec:group-algebras}
Let $G$ be a finite group.
Write $\IRR(G)$ for the set of irreducible complex characters of $G$, and set $d_\chi:=\chi(1_G)$ for $\chi\in\IRR(G)$.
The complex group algebra has the decomposition
\[
\C[G]\cong
\bigoplus_{\chi\in\IRR(G)}
M_{d_\chi}(\C)
\]
as associative $\C$-algebras \cite[p.~48, Proposition~10]{MR450380}.

For the real group algebra, we use the Frobenius--Schur theory recalled in \cite[p.~106]{MR450380}.
For each $\chi\in\IRR(G)$, define its Frobenius--Schur indicator by
\[
\nu(\chi):=\frac{1}{|G|}\sum_{g\in G}\chi(g^2).
\]
Then $\nu(\chi)\in\{1,0,-1\}$, and the corresponding irreducible complex representations are classified as follows.
\begin{enumerate}
\item[(F1)] If $\nu(\chi)=1$, then the representation is realizable over $\R$.
Such characters are said to be of real type.
\item[(F2)] If $\nu(\chi)=0$, then $\chi\neq\overline{\chi}$. 
Such characters are said to be of complex type.
\item[(F3)] If $\nu(\chi)=-1$, then $\chi=\overline{\chi}$, but the representation is not realizable over $\R$.
Such characters are said to be of quaternionic type.
In this case, $d_\chi$ is even.
\end{enumerate}
Here $\overline{\chi}$ denotes the character defined by $\overline{\chi}(g)=\overline{\chi(g)}$ for $g\in G$.

In the following decomposition, $\chi$ ranges over $\IRR(G)$.
For $\nu(\chi)=0$, we take exactly one character from each pair $\{\chi,\overline{\chi}\}$.
Then the real group algebra has the decomposition
\[
\R[G]\cong
\bigoplus_{\nu(\chi)=1} M_{d_\chi}(\R)
\oplus
\bigoplus_{\substack{\nu(\chi)=0\\\{\chi,\overline{\chi}\} }} M_{d_\chi}(\C)
\oplus
\bigoplus_{\nu(\chi)=-1} M_{d_\chi/2}(\BH)
\]
as associative $\R$-algebras, where $\BH$ denotes the real quaternion division algebra and the complex matrix algebras are regarded as real algebras.

The matrix size in the quaternionic case follows from restriction of scalars.
The corresponding irreducible real representation has real dimension $2d_\chi$ and commuting algebra $\BH$.
Its dimension over $\BH$ is therefore $d_\chi/2$, since $\dim_{\R}\BH=4$.

\section{Main results}
In this section, we first prove that adjacency algebras of Jordan schemes over fields of characteristic $0$ are semisimple.
We then determine the simple Jordan algebra decompositions over $\R$ or $\C$ for thin Jordan schemes.
We proceed in the order: non-regular, regular corresponding to a finite group, and regular corresponding to a non-associative RA-loop.

\subsection{Semisimplicity}
In this subsection, let $F$ be an arbitrary field of characteristic $0$.
It is known that the adjacency algebra of a Jordan scheme over $\R$ whose basis matrices are all symmetric is semisimple \cite{MR4462635}.
\begin{thm}\label{thm:semisimplicity}
Let $\JS$ be a Jordan scheme and let $\CA$ be its adjacency algebra over a field $\F$ of characteristic $0$.
Then $\CA$ is semisimple.
\end{thm}
\begin{proof}
We first show that the adjacency algebra $\R S$ is semisimple.
Assume that $\CI$ is a nonzero nil ideal of $\R S$, and choose $O \neq M \in \CI$.
Denote the $(i, j)$-entry of the matrix $M$ by $M_{i, j}$.
Since $\R S$ is closed under transposition, we have ${}^{t} M\in\R S$.
Since $\CI$ is an ideal, it follows that $M\star {}^{t} M\in\CI$.
For each $i$, we have $(M\star {}^{t} M)_{i, i} = \frac{1}{2}\sum_k\left(|M_{i, k}|^2+|M_{k, i}|^2\right) \geq 0$.
Since $M\neq O$, at least one entry of $M$ is nonzero.
Hence $\TRACE(M\star {}^{t} M) = \sum_{i, k}|M_{i, k}|^2 > 0$.
Therefore, $M\star {}^{t} M$ is not nilpotent as an element of $\R S$.
This contradicts the assumption that $\CI$ is a nil ideal.
Thus, $\R S$ is semisimple.

Let $F$ be an arbitrary field of characteristic $0$.
For $M\in FS$, let $T_M$ denote the product operator on $FS$.
Let $\tau_F$ denote the trace form on $FS$ from Subsection~\ref{subsec:jordan-algebra}, with the trace taken on the $(d+1)$-dimensional vector space $FS$.

By condition (JS3), the structure constants $p_{ij}^{k}$ lie in $\Q$ and are independent of the choice of $F$.
For each $k, r\in\{0, \ldots, d\}$, we have
\[
T_{A_k}(A_r)=\sum_{j=0}^{d}p_{kr}^{j}A_j.
\]
Thus, the diagonal entry indexed by $r$ in the matrix of $T_{A_k}$ with respect to the ordered basis $(A_0, \ldots, A_d)$ is $p_{kr}^{r}$.
Consequently, 
\[
\TRACE\left(T_{A_k}\right)=\sum_{r=0}^{d}p_{kr}^{r}.
\]
It follows that
\[
\tau_F(A_i, A_j) = \TRACE\left(T_{A_i\star A_j}\right) =\sum_{k=0}^{d}p_{ij}^{k}\TRACE\left(T_{A_k}\right) = \sum_{k=0}^{d}p_{ij}^{k}\sum_{r=0}^{d}p_{kr}^{r}.
\]
Therefore, the matrix of $\tau_F$ with respect to $(A_0, \ldots, A_d)$ is obtained from the rational matrix
\[
N:=\left(\sum_{k=0}^{d}p_{ij}^{k}\sum_{r=0}^{d}p_{kr}^{r}\right)_{0\leq i, j\leq d}\in M_{d+1}(\Q)
\]
by the natural embedding $\Q\hookrightarrow F$.

Since $\R S$ is semisimple, $\tau_{\R}$ is nondegenerate.
Hence $\det N\neq0$ over $\R$.
As $\det N$ is a rational number, it is a nonzero element of $\Q$.
Since $\F$ has characteristic $0$, its image under the embedding $\Q\hookrightarrow\F$ is also nonzero.
Thus, $\tau_{\F}$ is nondegenerate.
Applying the same trace-form criterion, we conclude that $\CA=\F S$ is semisimple.
\end{proof}

For the rest of this section, let $F$ be $\R$ or $\C$.
In what follows, we study the structure of the adjacency algebras of thin Jordan schemes.
We separate the cases: non-regular, regular corresponding to finite groups, and regular corresponding to non-associative RA-loops.

\subsection{Non-regular thin Jordan schemes}
We use the tensor-product description of $\CJ(\FG)$ established in Subsection~\ref{subsec:non-regular-thin-classification}.

\begin{thm}\label{thm:non-regular}
Let $\JS$ be a non-regular thin Jordan scheme, and let $\CA$ be its adjacency algebra over $\F$.
Let $G$ be a finite abelian group of order $n_G$, and let $P$ be a permutation matrix such that $\CA = P^{-1}\CJ(\FG)P$.
Let $n_{\RSYM}$ be the number of symmetric matrices in $S$.
Then
\[
\CA \cong 
\begin{cases}
\bigoplus_{i=1}^{n_G} \left (\SPIN_{\C}(2) \right), & \text{if $\F=\C$}, \\
\left(\bigoplus_{i=1}^{n_{\RSYM}} \left (\SPIN_{\R}(1, 1) \right)\right) \oplus 
\left(\bigoplus_{i=1}^{\frac{n_G-n_{\RSYM}}{2}}\left(\SPIN_{\C}(2) \right)_{\R} \right), & \text{if $\F=\R$}.
\end{cases}
\]
\end{thm}

\begin{proof}
By Proposition~\ref{prop:non-regular-classification}, there exists a permutation matrix $P$ such that $\CA=P^{-1}\CJ(\FG)P$.
We first show that $n_{\RSYM}=|G_2|$, where $G_2:=\{g\in G\mid g^2=1_G\}$.
The map $M\mapsto PMP^{-1}$, where $M\in\CA$, takes $S$ to the standard basis
\[
\left\{I_2\otimes P_g, \ E_{12}\otimes P_g, \ E_{21}\otimes P_g\mid g\in G\right\}
\]
of $\CJ(\FG)$.
Only matrices of the first form can be symmetric.
Since ${}^tP_g=P_{g^{-1}}$, the matrix $I_2\otimes P_g$ is symmetric if and only if $g^2=1_G$.
As permutation conjugation preserves symmetry, we obtain $n_{\RSYM}=|G_2|$.

Recall that
$
\CA\cong\CJ(\FG)\cong \CJF \otimes_{\F}\FG.
$
Set $D_1:=E_{12}+E_{21}$ and $D_2:=E_{12}-E_{21}$.
Then $\{I_2, D_1, D_2\}$ is an $\F$-basis of $\CJF$, with identity element $I_2$, and $D_1^{\star 2}=I_2, \quad D_2^{\star 2}=-I_2, \quad D_1\star D_2=O$.

\textit{The case $\F=\C$.}
Since $G$ is abelian, all its irreducible complex characters are linear.
For each $\chi\in\IRR(G)$, define
\[
e_\chi:=\frac{1}{|G|}\sum_{g\in G}\overline{\chi(g)}g.
\]
These elements form a complete set of pairwise orthogonal primitive idempotents of $\C[G]$ \cite[Chapter~3, \S~2]{MR998775}.
Since $ge_\chi=\chi(g)e_\chi$ for every $g\in G$, we have
$
\C[G]e_\chi=\C e_\chi=\{ze_\chi\mid z\in\C\}.
$
Hence
\[
\C[G]=\bigoplus_{\chi\in\IRR(G)}\C e_\chi\cong\bigoplus_{i=1}^{n_G}\C.
\]
Since $e_\chi^2=e_\chi$, the map 
\[
\CJ(\C)\otimes_{\C}\C e_\chi\longrightarrow \CJ(\C), \quad M\otimes(ze_\chi)\longmapsto zM
\]
is an isomorphism of complex Jordan algebras.

The basis $\{I_2, D_1, \RI D_2\}$ of $\CJ(\C)$ satisfies $D_1^{\star 2}=(\RI D_2)^{\star 2}=I_2$ and $D_1\star(\RI D_2)=O$.
Hence 
\[
\CJ(\C) \cong \SPIN_{\C}(2), 
\]
and consequently
\[
\CA \cong \CJ(\C)\otimes_{\C}\C[G] \cong\bigoplus_{\chi\in\IRR(G)} \left(\CJ(\C)\otimes_{\C}\C e_\chi\right) \cong\bigoplus_{i=1}^{n_G}\CJ(\C) \cong\bigoplus_{i=1}^{n_G}\left(\SPIN_{\C} (2) \right).
\]

\textit{The case $\F=\R$.}
Let $\widehat{G}$ be the character group of the finite abelian group $G$.
Since $\widehat{G}\cong G$ and $\chi=\overline{\chi}$ if and only if $\chi^2$ is the trivial character, there are $|G_2|=n_{\RSYM}$ self-conjugate characters.
The remaining $n_G-n_{\RSYM}$ characters form $(n_G-n_{\RSYM})/2$ conjugate pairs.

If $\chi=\overline{\chi}$, then $e_\chi\in\R[G]$ and
$
\R[G]e_\chi=\R e_\chi=\{re_\chi\mid r\in\R\}\cong\R.
$
Thus, a self-conjugate character gives a one-dimensional real ideal.

If $\chi\neq\overline{\chi}$, set $e_{[\chi]}:=e_\chi+e_{\overline{\chi}}$.
Since complex conjugation sends $e_\chi$ to $e_{\overline{\chi}}$, we have $e_{[\chi]}\in\R[G]$ and
$
\R[G]e_{[\chi]}=\{ze_\chi+\overline{z}e_{\overline{\chi}}\mid z\in\C\}\cong\C, 
$
where the isomorphism is given by
$ze_\chi+\overline{z}e_{\overline{\chi}}\mapsto z$.
Here $\C$ is regarded as a real algebra.
In particular, this ideal is two-dimensional over $\R$.
Therefore, 
\[
\R[G]=\left(\bigoplus_{\chi=\overline{\chi}}\R e_\chi\right)\oplus\left(\bigoplus_{\chi\neq\overline{\chi}}\R[G]e_{[\chi]}\right), 
\]
where the first sum runs through all self-conjugate irreducible characters and the second through one representative from each pair of non-self-conjugate irreducible characters.
Thus, as real algebras, 
\[
\R[G]\cong\left(\bigoplus_{i=1}^{n_{\RSYM}}\R\right)\oplus\left(\bigoplus_{i=1}^{\frac{n_G-n_{\RSYM}}{2}}\C\right).
\]

For a self-conjugate character $\chi$, the map $M\otimes(re_\chi)\mapsto rM$ gives an isomorphism $\CJ(\R)\otimes_{\R}\R e_\chi\cong \CJ(\R)$.
The product relations for the basis $\{I_2, D_1, D_2\}$ give
\[
\CJ(\R) \cong \SPIN_{\R} (1, 1) .
\]

The map $M\otimes z\mapsto zM$ is an isomorphism of real Jordan algebras from $\CJ(\R)\otimes_{\R}\C$ onto the realification of the complex Jordan algebra $\CJ(\C)$.
By the complex case, we therefore have an isomorphism of real Jordan algebras
\[
\CJ(\R) \otimes_{\R} \C \cong \left( \SPIN_{\C} (2) \right)_{\R}.
\]
Consequently, 
\[
\begin{aligned}
\CA
&\cong \CJ(\R) \otimes_{\R}\R[G]\\
&\cong \left(\bigoplus_{i=1}^{n_{\RSYM}} \CJ(\R) \right) \oplus \left(\bigoplus_{i=1}^{\frac{n_G-n_{\RSYM}}{2}} \left( \CJ(\R) \otimes_{\R} \C \right) \right)\\
&\cong \left( \bigoplus_{i=1}^{n_{\RSYM}} \left( \SPIN_{\R} (1, 1) \right) \right) \oplus \left( \bigoplus_{i=1}^{\frac{n_G-n_{\RSYM}}{2}} \left(\SPIN_{\C} (2) \right)_{\R}\right).
\end{aligned}
\]
\end{proof}

\subsection{Regular thin Jordan schemes corresponding to finite groups}
We determine the adjacency algebras of regular thin Jordan schemes corresponding to finite groups over $\R$ or $\C$.
We use the group-algebra decompositions recalled in Subsection~\ref{subsec:group-algebras}.

For an associative algebra $\CB$ over $\F$, write $\CB^{+}$ for the Jordan algebra with the same underlying vector space and product $x\star y=(xy+yx)/2$.

Let $G$ be a finite group.
For each $g\in G$, let $P_g$ denote the permutation matrix of the left translation $x\mapsto gx$, with the convention $P_g\bm{e}_x=\bm{e}_{gx}\quad(x\in G)$, where $\{\bm{e}_x\mid x\in G\}$ is the standard basis of $\F^G$.
Write $\CA(G,G):=\SPAN_{\F}\{P_g\mid g\in G\}$ for the adjacency algebra of the thin association scheme associated with $G$, equipped with ordinary matrix product.
Since $P_gP_h=P_{gh}$ for all $g,h\in G$, the linear map
\[
\Phi_\F\colon\FG\longrightarrow\CA(G,G),\quad\sum_{g\in G}a_gg\longmapsto\sum_{g\in G}a_gP_g
\]
preserves ordinary product.
It is surjective by the definition of $\CA(G,G)$.
Moreover, the matrices $P_g$ are linearly independent, since
\[
\left(\sum_{g\in G}a_gP_g\right)\bm{e}_{1_G}=\sum_{g\in G}a_g\bm{e}_g.
\]
Thus, $\Phi_\F$ is an associative algebra isomorphism.
The same linear map therefore induces a Jordan algebra isomorphism
\[
\Phi_\F\colon\FG^{+}\longrightarrow\CA(G,G)^{+}.
\]

Let $\JS$ be a regular thin Jordan scheme corresponding to $G$, and let $\CA$ be its adjacency algebra over $\F$.
After a suitable relabelling of the points of $X$, its adjacency matrices are precisely the matrices $\{P_g\mid g\in G\}$.
Consequently, $\CA\cong\CA(G,G)^{+}\cong\FG^{+}$ as Jordan algebras.
Combining this identification with the group-algebra decompositions recalled in Subsection~\ref{subsec:group-algebras} yields the following theorem.

\begin{thm}\label{thm:regular-thin-group}
Let $\JS$ be a regular thin Jordan scheme corresponding to a finite group $G$, and let $\CA$ be its adjacency algebra over $\F$.

If $\F = \C$, then
\[
\CA \cong \bigoplus_{\chi\in\IRR(G)} M_{d_\chi}(\C)^+.
\]

If $\F = \R$, then
\[
\CA \cong 
\bigoplus_{\nu(\chi)=1} M_{d_\chi}(\R)^+
\oplus
\bigoplus_{\substack{\nu(\chi)=0\\\{\chi,\overline{\chi}\} }} \bigl(M_{d_\chi}(\C)^+\bigr)_{\R}
\oplus
\bigoplus_{\nu(\chi)=-1} M_{d_\chi/2}(\BH)^+.
\]
In the real decomposition, $\chi$ ranges over $\IRR(G)$, with exactly one character chosen from each pair $\{\chi,\overline{\chi}\}$ when $\nu(\chi)=0$.
Every summand in these decompositions is simple over the indicated ground field.
\end{thm}

\begin{proof}
Let $\CB$ be a finite-dimensional semisimple associative algebra over $\F$, where $\CHAR(\F)\neq2$, and let
$
\CB=\bigoplus_{i=1}^{r}\CB_i
$
be its decomposition into simple associative ideals.
For distinct indices $i$ and $j$, we have $\CB_i\CB_j=\CB_j\CB_i=0$.
Consequently, each $\CB_i$ is also an ideal under the Jordan product, and
$
\CB^{+}=\bigoplus_{i=1}^{r}\CB_i^{+}.
$
By \cite[p.~169, Full Example~3.1.1(6)]{MR2014924}, simplicity of each associative algebra $\CB_i$ implies simplicity of the Jordan algebra $\CB_i^{+}$.

Applying this observation to the associative decompositions of $\C[G]$ and $\R[G]$ gives the asserted decompositions of $\C[G]^{+}$ and $\R[G]^{+}$.
For the complex matrix components regarded as real algebras, restriction of scalars commutes with the passage to the Jordan product:
\[
\bigl(M_{d_\chi}(\C)_{\R}\bigr)^{+}
=
\bigl(M_{d_\chi}(\C)^{+}\bigr)_{\R}.
\]
Every resulting summand is therefore simple over the indicated ground field.
\end{proof}

In the real decomposition, the dimensions of the summands corresponding to real, complex, and quaternionic types are, respectively,
\[
d_\chi^2, \quad 2d_\chi^2, \quad 4\left(\frac{d_\chi}{2}\right)^2=d_\chi^2.
\]

\subsection{Regular thin Jordan schemes corresponding to non-associative RA-loops}
In this subsection, we determine the adjacency algebras of regular thin Jordan schemes corresponding to non-associative RA-loops over $\R$ or $\C$ under Theorem~\ref{thm:regular-classification}.
Let $\JS$ be such a scheme, and let $\CA$ be its adjacency algebra over $\F$.
Choose a finite group $G$ and an element $g_0\in Z(G)$ such that the corresponding non-associative RA-loop is isomorphic to $L(G, *, g_0)$.
In what follows, we write $Z := Z(G)$ and $L := L(G, *, g_0)$.
Let $n_G := |G|$, $n_Z := |Z| = n_G/4$, and $Z = \{z_1, z_2, \ldots, z_{n_Z}\}$.
Let $G' = \{1_G, s\}$ be the commutator subgroup of $G$, where $s$ is a central involution.
Choose $a, b\in G$ such that
$
G/Z = \{Z, aZ, bZ, abZ\}.
$

For each $z\in Z$, let $P_z\in M_{n_Z}(\R)$ be the permutation matrix corresponding to left product by $z$ on $Z$, with respect to the ordering $(z_1, z_2, \ldots, z_{n_Z})$.
More precisely, the $(i, j)$-entry of $P_z$ is given by
\[
(P_z)_{ij} =
\begin{cases}
1 & \text{if } z_i = zz_j, \\
0 & \text{otherwise}.
\end{cases}
\]

For each $c\in\{1_G, a, b, ab\}$ and $x\in\{0, 1\}$, order the elements of $cZ\times\{x\}$ as
\[
\bigl((cz_1, x), (cz_2, x), \ldots, (cz_{n_Z}, x)\bigr).
\]
We then order the elements of $L$ as
\[
\bigl( Z\times\{0\}, aZ\times\{0\}, bZ\times\{0\}, abZ\times\{0\}, Z\times\{1\}, aZ\times\{1\}, bZ\times\{1\}, abZ\times\{1\} \bigr).
\]
With respect to this ordering, the left-translation matrices admit the following reduction.

\begin{lemm}\label{lem:left-translations}
For $z\in Z$, we have $\ell_{(z, 0)} = I_8 \otimes P_z$, where $\otimes$ denotes the Kronecker product of matrices.
Then $\ell_{(cz, x)} = \ell_{(z, 0)}\ell_{(c, x)} \quad (c\in\{1_G, a, b, ab\}, \ x\in\{0, 1\})$.
Consequently, all the left-translation matrices are determined by the matrices 
$\ell_{(z, 0)}$ ($z\in Z$) and the seven matrices 
\[
\ell_{(a, 0)}, \ \ell_{(b, 0)}, \ \ell_{(ab, 0)}, \ \ell_{(1_G, 1)}, \ \ell_{(a, 1)}, \ \ell_{(b, 1)}, \ \ell_{(ab, 1)}.
\]
Moreover, the matrices $\{\ell_u \mid u \in L\}$ form a basis of $\CA$.

Their block-matrix forms are given below.
Each entry in the following displays is an $n_Z \times n_Z$ matrix, and each blank entry denotes the $n_Z \times n_Z$ zero matrix.
In particular, $I_{n_Z}$ denotes the $n_Z \times n_Z$ identity matrix.
First, 
{\fontsize{9pt}{11pt}\selectfont\begin{equation*}
\setlength{\arraycolsep}{3.5pt}
\renewcommand{\arraystretch}{1.1}
\displaystyle
\ell_{(z, 0)}
 =
\left[
\begin{array}{cccc|cccc}
P_z& & & & & & & \\
& P_z& & & & & & \\
& & P_z& & & & & \\
& & & P_z& & & & \\ \hline
& & & & P_z& & & \\
& & & & & P_z& & \\
& & & & & & P_z& \\
& & & & & & & P_z
\end{array}
\right], 
\displaystyle
\ell_{(1_G, 1)}
 = 
\left[
\begin{array}{cccc|cccc}
& & & & P_{g_0}& & & \\
& & & & & P_{g_0s}& & \\
& & & & & & P_{g_0s}& \\
& & & & & & & P_{g_0s}\\ \hline
I_{n_Z}& & & & & & & \\
& P_s& & & & & & \\
& & P_s& & & & & \\
& & & P_s& & & & 
\end{array}
\right].
\end{equation*}}
The matrices associated with $a$ are
{\fontsize{9pt}{11pt}\selectfont\begin{equation*}
\setlength{\arraycolsep}{2pt}
\renewcommand{\arraystretch}{1.1}
\displaystyle
\ell_{(a, 0)}
 = 
\left[
\begin{array}{cccc|cccc}
& P_{a^2}& & & & & & \\
I_{n_Z}& & & & & & & \\
& & & P_{a^2}& & & & \\
& & I_{n_Z}& & & & & \\ \hline
& & & & & P_{a^2}& & \\
& & & & I_{n_Z}& & & \\
& & & & & & & P_{a^2s}\\
& & & & & & P_s& 
\end{array}
\right], 
\displaystyle
\ell_{(a, 1)}
 = 
\left[
\begin{array}{cccc|cccc}
& & & & & P_{g_0a^2s}& & \\
& & & & P_{g_0}& & & \\
& & & & & & & P_{g_0a^2}\\
& & & & & & P_{g_0}& \\ \hline
& P_{a^2s}& & & & & & \\
I_{n_Z}& & & & & & & \\
& & & P_{a^2s}& & & & \\
& & P_s& & & & & 
\end{array}
\right].
\end{equation*}}
The matrices associated with $b$ are
{\fontsize{9pt}{11pt}\selectfont\begin{equation*}
\setlength{\arraycolsep}{2pt}
\renewcommand{\arraystretch}{1.1}
\displaystyle
\ell_{(b, 0)}
 = 
\left[
\begin{array}{cccc|cccc}
& & P_{b^2}& & & & & \\
& & & P_{b^2s}& & & & \\
I_{n_Z}& & & & & & & \\
& P_s& & & & & & \\ \hline
& & & & & & P_{b^2}& \\
& & & & & & & P_{b^2}\\
& & & & I_{n_Z}& & & \\
& & & & & I_{n_Z}& & 
\end{array}
\right], 
\displaystyle
\ell_{(b, 1)}
 = 
\left[
\begin{array}{cccc|cccc}
& & & & & & P_{g_0b^2s}& \\
& & & & & & & P_{g_0b^2s}\\
& & & & P_{g_0}& & & \\
& & & & & P_{g_0s}& & \\ \hline
& & P_{b^2s}& & & & & \\
& & & P_{b^2}& & & & \\
I_{n_Z}& & & & & & & \\
& I_{n_Z}& & & & & & 
\end{array}
\right].
\end{equation*}}
Finally, the matrices associated with $ab$ are
{\fontsize{9pt}{11pt}\selectfont\begin{equation*}
\setlength{\arraycolsep}{0pt}
\renewcommand{\arraystretch}{1.1}
\displaystyle
\ell_{(ab, 0)} = 
\left[
\begin{array}{cccc|cccc}
& & & P_{a^2b^2s}& & & & \\
& & P_{b^2}& & & & & \\
& P_{a^2s}& & & & & & \\
I_{n_Z}& & & & & & & \\ \hline
& & & & & & & P_{a^2b^2s}\\
& & & & & & P_{b^2s}& \\
& & & & & P_{a^2}& & \\
& & & & I_{n_Z}& & & 
\end{array}
\right], 
\ell_{(ab, 1)} = 
\left[
\begin{array}{cccc|cccc}
& & & & & & & P_{g_0a^2b^2}\\
& & & & & & P_{g_0b^2}& \\
& & & & & P_{g_0a^2s}& & \\
& & & & P_{g_0}& & & \\ \hline
& & & P_{a^2b^2}& & & & \\
& & P_{b^2s}& & & & & \\
& P_{a^2}& & & & & & \\
I_{n_Z}& & & & & & & 
\end{array}
\right].
\end{equation*}}
\end{lemm}
\begin{proof}
We illustrate the block-matrix calculation by determining one block of $\ell_{(az, 1)}$, where $z\in Z$.
For $w\in Z$, using $(abw)^*=abws$, $ba=abs$, and $s^2=1_G$, we obtain
\[
\begin{aligned}
\lambda_{(az, 1)}(abw, 1)
&= (az, 1)\diamond(abw, 1)\\
&= (g_0(abw)^*az, 0)\\
&= (g_0abwsaz, 0)\\
&= (g_0a^2bws^2z, 0)\\
&= (b(g_0a^2z)w, 0).
\end{aligned}
\]
Thus, the entry of $\ell_{(az, 1)}$ in the row indexed by $(b(g_0a^2z)w, 0)$ and the column indexed by $(abw, 1)$ is $1$.
As $w$ ranges over $Z$, the block with rows indexed by $bZ\times\{0\}$ and columns indexed by $abZ\times\{1\}$ represents the permutation $w\mapsto(g_0a^2z)w$ of $Z$.
Hence this block is $P_{g_0a^2z}$.
In particular, setting $z=1_G$ gives the $(3, 8)$-block $P_{g_0a^2}$ of $\ell_{(a, 1)}$ displayed above.
The other blocks are obtained similarly.

A direct matrix calculation gives
\[
\ell_{(z, 0)}\ell_{(c, x)} = \ell_{(cz, x)} \quad (z\in Z, \ c\in\{1_G, a, b, ab\}, \ x\in\{0, 1\}).
\]
With respect to the ordering above, we also have $\ell_{(z, 0)} = I_8\otimes P_z$.

By \cite[Theorem~3.17]{MR5094503}, the matrices $\{\ell_u\mid u\in L\}$ are the adjacency matrices of a thin Jordan scheme.
Thus, they are nonzero $\{0, 1\}$-matrices with pairwise disjoint supports, and hence they are linearly independent.
Since they span $\CA$ by definition, they form a basis of $\CA$.
\end{proof}

Define maps $\rho\colon L\to\F L$ and $\sigma\colon Z\to\F L$ by
\[
\rho((g, x)) := \frac{1}{2}\bigl((g, x)+(gs, x)\bigr), \quad \sigma(z) := \frac{1}{2}\bigl((z, 0)-(zs, 0)\bigr)
\]
for $(g, x)\in L$ and $z\in Z$.
Set $H := \{\ell_{\rho((g, x))}\mid(g, x)\in L\}, K := \{\ell_{\sigma(z)}\mid z\in Z\}$.

\begin{lemm}\label{lem:central-elements-regular}
For $(g_1, x_1), (g_2, x_2)\in L$ and $z\in Z$, the following identities hold:
\[
\ell_{\rho((g_1, x_1))}, \, \ell_{\sigma(z)}\in Z(\CA), 
\]
\[
\ell_{(z, 0)}\ell_{(g_1, x_1)} = \ell_{(g_1z, x_1)} = \ell_{(g_1, x_1)}\ell_{(z, 0)}, 
\]
\[
\ell_{\rho((g_1, x_1))}\ell_{(g_2, x_2)} = \ell_{(g_2, x_2)}\ell_{\rho((g_1, x_1))} = \ell_{\rho((g_1, x_1))}\ell_{\rho((g_2, x_2))} = \ell_{\rho((g_2, x_2)\diamond(g_1, x_1))}, 
\]
\[
\ell_{\rho((g_1, x_1))}\star \ell_{\sigma(z)} = O.
\]
\end{lemm}
\begin{proof}
We use the sufficient condition for membership in the center given in Subsection~\ref{subsec:jordan-algebra}.
A direct matrix calculation shows that, for $z\in Z$ and $(g_1, x_1)\in L$, $\ell_{(z, 0)}\ell_{(g_1, x_1)} = \ell_{(g_1z, x_1)} = \ell_{(g_1, x_1)}\ell_{(z, 0)}$.
Thus, $\ell_{(z, 0)}\in Z(\CA)$.
Since $\ell_{\sigma(z)} = \frac{1}{2}\left(\ell_{(z, 0)}-\ell_{(zs, 0)}\right)$, we also have $\ell_{\sigma(z)}\in Z(\CA)$.

Next, let $(g_1, x_1), (g_2, x_2), (g_3, x_3)\in L$.
Substituting the four cases in the definition of $\diamond$ shows that each of the operators 
$\ell_{\rho((g_1, x_1))}\ell_{(g_2, x_2)}$, $\ell_{(g_2, x_2)}\ell_{\rho((g_1, x_1))}$, 
$\ell_{\rho((g_1, x_1))}\ell_{\rho((g_2, x_2))}$, 
$\ell_{\rho((g_2, x_2)\diamond(g_1, x_1))}$ 
maps the basis element $(g_3, x_3)$ of $\F L$ to
\[
\rho\left((g_0^{\left\lfloor(x_1+x_2+x_3)/2\right\rfloor}g_1g_2g_3, \, x_1+x_2+x_3\pmod{2})\right).
\]
Hence $\ell_{\rho((g_1, x_1))}\ell_{(g_2, x_2)} = \ell_{(g_2, x_2)}\ell_{\rho((g_1, x_1))} = \ell_{\rho((g_1, x_1))}\ell_{\rho((g_2, x_2))} = \ell_{\rho((g_2, x_2)\diamond(g_1, x_1))}$.
Since the matrices $\ell_{(g_2, x_2)}$ span $\CA$, it follows that $\ell_{\rho((g_1, x_1))}\in Z(\CA)$.

Finally, using the preceding commutativity, we obtain
\[
\begin{aligned}
4\left(\ell_{\rho((g_1, x_1))}\star \ell_{\sigma(z)}\right)& = \left(\ell_{(g_1, x_1)}+\ell_{(g_1s, x_1)}\right)\left(\ell_{(z, 0)}-\ell_{(zs, 0)}\right)\\& = \ell_{(g_1z, x_1)}-\ell_{(g_1zs, x_1)}+\ell_{(g_1zs, x_1)}-\ell_{(g_1z, x_1)} = O.
\end{aligned}
\]
Therefore, $\ell_{\rho((g_1, x_1))}\star \ell_{\sigma(z)} = O$.
\end{proof}

\begin{lemm}\label{lem:HK}
The sets $(H, \star)$ and $(K, \star)$ are abelian groups.
Their identity elements are $\ell_{\rho((1_G, 0))}$ and $\ell_{\sigma(1_G)}$, respectively.
Moreover, $(K, \star)\cong Z$.
\end{lemm}
\begin{proof}
By Lemma~\ref{lem:central-elements-regular}, the product in $H$ is associative and commutative, and satisfies $\ell_{\rho((g, x))}\star \ell_{\rho((h, y))} = \ell_{\rho((h, y)\diamond(g, x))}$.
The identity element of $H$ is $\ell_{\rho((1_G, 0))}$, and the inverse of $\ell_{\rho((g, x))}$ is $\ell_{\rho((g^{-1}g_0^{-x}, x))}$.
For $z, w\in Z$, we have $\ell_{\sigma(z)}\star \ell_{\sigma(w)} = \ell_{\sigma(zw)}$.
Thus, the identity element of $K$ is $\ell_{\sigma(1_G)}$, and the inverse of $\ell_{\sigma(z)}$ is $\ell_{\sigma(z^{-1})}$.
Consequently, the map $Z\longrightarrow K, z\longmapsto \ell_{\sigma(z)}$, is a group isomorphism.
\end{proof}

Via the isomorphism $Z\cong K$ in Lemma~\ref{lem:HK}, we use characters of $Z$ to describe the idempotents associated with $K$.
For $\phi\in\IRR(H)$ and $\psi\in\IRR(Z)$, define
\[
e_\phi := \frac{1}{|H|}\sum_{h\in H}\overline{\phi(h)}h, \quad f_\psi := \frac{1}{|Z|}\sum_{z\in Z}\overline{\psi(z)}\ell_{\sigma(z)}.
\]

\begin{lemm}\label{lem:nonzero-idempotents-regular}
For every $\phi\in\IRR(H)$, we have $e_\phi\neq O$.
Moreover, for $\psi\in\IRR(Z)$, we have $f_\psi\neq O$ if and only if $\psi(s) = -1$.
\end{lemm}

\begin{proof}
Since the elements of $H$ are linearly independent in $\CA$, we have $e_\phi\neq O$ for every $\phi\in\IRR(H)$.
Let $\{w_1, w_2, \ldots, w_{n_Z/2}\}$ be a set of representatives of $Z/\langle s\rangle$.
Since $\ell_{\sigma(w_i s)} = -\ell_{\sigma(w_i)}$, we have
\[
f_\psi 
= \frac{1}{|Z|}\sum_{i = 1}^{n_Z/2}\left(\overline{\psi(w_i)}\ell_{\sigma(w_i)}+\overline{\psi(w_i s)}\ell_{\sigma(w_i s)}\right) 
= \frac{1}{|Z|}\sum_{i = 1}^{n_Z/2}\left(1-\overline{\psi(s)}\right)\overline{\psi(w_i)}\ell_{\sigma(w_i)}.
\]
Since $s^2 = 1_G$, we have $\psi(s)\in\{1, -1\}$.
Thus, $f_\psi = O$ if $\psi(s) = 1$.
If $\psi(s) = -1$, the linear independence of $\ell_{\sigma(w_1)}, \ldots, \ell_{\sigma(w_{n_Z/2})}$ implies $f_\psi\neq O$.
\end{proof}

When $\F = \R$, regard $e_\phi$ and $f_\psi$ as elements of the complexification of $\CA$, and define
\[
e_{[\phi]} := \begin{cases} e_\phi, & \text{if $\phi = \overline{\phi}$}, \\e_\phi+e_{\overline{\phi}}, & \text{if $\phi\neq\overline{\phi}$}, \end{cases} \quad 
f_{[\psi]} := \begin{cases} f_\psi, & \text{if $\psi = \overline{\psi}$}, \\f_\psi+f_{\overline{\psi}}, & \text{if $\psi\neq\overline{\psi}$}.\end{cases}
\]

In all sums and direct sums below for $\F = \R$, $\phi\in\IRR(H)$ is understood to range over all self-conjugate characters and exactly one character from each pair $\{\phi, \overline{\phi}\}$ of non-self-conjugate characters.
Similarly, $\psi\in\IRR(Z)$ with $\psi(s)\neq1$ is understood to range over all self-conjugate characters satisfying this condition and exactly one character from each pair $\{\psi, \overline{\psi}\}$ of non-self-conjugate characters satisfying this condition.

\begin{prop}\label{prop:central-idempotents-regular}
For the chosen field $\F$, the corresponding elements defined above are pairwise orthogonal central idempotents of $\CA$ whose sum is the identity element $\ell_{(1_G, 0)}$.
Consequently, if $\F = \C$, then
\[
\CA = \left(\bigoplus_{\phi\in\IRR(H)}\left(\CA \star e_\phi\right)\right)\oplus\left(\bigoplus_{\substack{\psi\in\IRR(Z)\\\psi(s)\neq1}}\left(\CA \star f_\psi\right)\right), 
\]
whereas if $\F = \R$, then
\[
\CA = \left(\bigoplus_{\phi\in\IRR(H)}\left(\CA \star e_{[\phi]}\right)\right)\oplus\left(\bigoplus_{\substack{\psi\in\IRR(Z)\\\psi(s)\neq1}}\left(\CA \star f_{[\psi]}\right)\right).
\]
\end{prop}

\begin{proof}
By Lemma~\ref{lem:central-elements-regular}, $H, K \subseteq Z(\CA)$.
Hence all the elements defined above are central.
Their idempotency and mutual orthogonality follow from the character theory of the abelian groups $H$ and $Z$, the isomorphism $Z\cong K$, and $\ell_{\rho((g, x))}\star \ell_{\sigma(z)} = O$.

If $\F = \C$, character orthogonality gives
\[
\sum_{\phi\in\IRR(H)}e_\phi + \sum_{\substack{\psi\in\IRR(Z)\\\psi(s)\neq1}}f_\psi = \frac{1}{2}\left(\ell_{(1_G, 0)}+\ell_{(s, 0)}\right)+\frac{1}{2}\left(\ell_{(1_G, 0)}-\ell_{(s, 0)}\right) = \ell_{(1_G, 0)}.
\]
Thus, the idempotents sum to the identity element of $\CA$.

If $\F = \R$, the same conclusion follows by combining each pair of complex-conjugate idempotents.
The stated direct-sum decompositions now follow.
\end{proof}

Finally, we determine the Jordan algebra structure of each direct summand.

\begin{thm}\label{main:regular}
Let $\JS$ be a regular thin Jordan scheme corresponding to a non-associative RA-loop and let $\CA$ be its adjacency algebra over $\F$.
Choose a finite group $G$ and an element $g_0\in Z(G)$ such that the corresponding RA-loop is isomorphic to $L(G, *, g_0)$.
Let $H$ be the abelian group defined above.
Then the following statements hold.

If $\F = \C$, then
\[
\CA\cong\left(\bigoplus_{\phi\in\IRR(H)}\C\right)
\oplus 
\left(\bigoplus_{\substack{\psi\in\IRR(Z)\\\psi(s)\neq1}} \SPIN_{\C} (7) \right).
\]

If $\F=\R$, then
\[
\begin{aligned}
\CA \cong{}
& \left( \bigoplus_{\substack{\phi\in\IRR(H)\\ \phi=\overline{\phi}}} \R \right)
\oplus
\left( \bigoplus_{\substack{\phi\in\IRR(H)\\ \phi\neq\overline{\phi}}} \C \right)
\oplus
\left( \bigoplus_{\substack{\psi\in\IRR(Z)\\ \psi(s)\neq1, \ \psi\neq\overline{\psi}}} \left(\SPIN_{\C} (7) \right)_{\R} \right)\\
&\oplus
\left( \bigoplus_{\substack{\psi\in\IRR(Z)\\ \psi(s)\neq1, \ \psi=\overline{\psi}\\ \psi(a^2)=\psi(b^2)=\psi(g_0)=-1}} \SPIN_{\R} (0, 7) \right)
\oplus
\left( \bigoplus_{\substack{\psi\in\IRR(Z)\\ \psi(s)\neq1, \ \psi=\overline{\psi}\\ (\psi(a^2), \psi(b^2), \psi(g_0))\neq(-1, -1, -1)}} \SPIN_{\R} (4, 3) \right).
\end{aligned}
\]
\end{thm}
\begin{proof}
We treat the complex and real cases separately.
Define elements of $L$ as 
\[
q_1 := (a, 0), q_2 := (b, 0), q_3 := (ab, 0), 
\]
\[
q_4 := (1_G, 1), q_5 := (a, 1), q_6 := (b, 1), q_7 := (ab, 1).
\]
For each $\psi\in\IRR(Z)$ with $\psi(s)\neq1$, set $\alpha := \psi(a^2), \beta := \psi(b^2), \gamma := \psi(g_0)$, noting that $a^2, b^2, g_0\in Z$.

\textit{The case $\F = \C$.}
For $(g, x)\in L$, set $h_{g, x} := \ell_{\rho((g, x))}\in H$.
By Lemma~\ref{lem:central-elements-regular}, we have $\ell_{(g, x)}\star e_\phi = h_{g, x}\star e_\phi$.
Hence $\ell_{(g, x)}\star e_\phi = \phi(h_{g, x})e_\phi$.
Since the matrices $\ell_{(g, x)}$ form a basis of $\CA$, it follows that
\[
\CA \star e_\phi = \C e_\phi\cong\C.
\]

Let $\psi\in\IRR(Z)$ satisfy $\psi(s)\neq1$.
Then $\psi(s) = -1$.
For $1\leq i\leq7$, set $v_i := \ell_{q_i}\star f_\psi$.
To write these elements explicitly, define the $n_Z\times n_Z$ matrix
\[
Q_\psi := \frac{1}{2n_Z}\sum_{z\in Z}\overline{\psi(z)}(P_z-P_{zs}) = \frac{1}{n_Z}\sum_{z\in Z}\overline{\psi(z)}P_z, 
\]
where the second equality follows from $\psi(s)=-1$.
For every $w\in Z$, reindexing the sum gives
\[
P_wQ_\psi = \frac{1}{n_Z}\sum_{z\in Z}\overline{\psi(w^{-1}z)}P_z = \psi(w)Q_\psi.
\]
The matrices $P_z$ are linearly independent, so $Q_\psi\neq O$.
The preceding identity gives $Q_\psi^2=Q_\psi$, and $f_\psi=I_8\otimes Q_\psi$ follows from its definition.
Moreover, $f_\psi$ commutes with each $\ell_{q_i}$ under the ordinary matrix product, so $v_i=\ell_{q_i}f_\psi$.
By Lemma~\ref{lem:left-translations}, the eight matrices have the following block forms.
In each $8\times8$ scalar matrix, blank entries denote $0$.
\begin{equation*}
\setlength{\arraycolsep}{2pt}
\renewcommand{\arraystretch}{1.1}
f_\psi =
\left[
\begin{array}{cccc|cccc}
1& & & & & & & \\
& 1& & & & & & \\
& & 1& & & & & \\
& & & 1& & & & \\ \hline
& & & & 1& & & \\
& & & & & 1& & \\
& & & & & & 1& \\
& & & & & & & 1
\end{array}
\right]\otimes Q_\psi, \quad
v_1 =
\left[
\begin{array}{cccc|cccc}
& \alpha& & & & & & \\
1& & & & & & & \\
& & & \alpha& & & & \\
& & 1& & & & & \\ \hline
& & & & & \alpha& & \\
& & & & 1& & & \\
& & & & & & & -\alpha \\
& & & & & & -1& 
\end{array}
\right]\otimes Q_\psi, 
\end{equation*}
\begin{equation*}
\setlength{\arraycolsep}{2pt}
\renewcommand{\arraystretch}{1.1}
v_2 =
\left[
\begin{array}{cccc|cccc}
& & \beta& & & & & \\
& & & -\beta& & & & \\
1& & & & & & & \\
& -1& & & & & & \\ \hline
& & & & & & \beta& \\
& & & & & & & \beta \\
& & & & 1& & & \\
& & & & & 1& & 
\end{array}
\right]\otimes Q_\psi, \quad
v_3 =
\left[
\begin{array}{cccc|cccc}
& & & -\alpha\beta& & & & \\
& & \beta& & & & & \\
& -\alpha& & & & & & \\
1& & & & & & & \\ \hline
& & & & & & & -\alpha\beta \\
& & & & & & -\beta& \\
& & & & & \alpha& & \\
& & & & 1& & & 
\end{array}
\right]\otimes Q_\psi, 
\end{equation*}
\begin{equation*}
\setlength{\arraycolsep}{1pt}
\renewcommand{\arraystretch}{1.1}
v_4 =
\left[
\begin{array}{cccc|cccc}
& & & & \gamma& & & \\
& & & & & -\gamma& & \\
& & & & & & -\gamma& \\
& & & & & & & -\gamma \\ \hline
1& & & & & & & \\
& -1& & & & & & \\
& & -1& & & & & \\
& & & -1& & & & 
\end{array}
\right]\otimes Q_\psi, \quad
v_5 =
\left[
\begin{array}{cccc|cccc}
& & & & & -\gamma\alpha& & \\
& & & & \gamma& & & \\
& & & & & & & \gamma\alpha \\
& & & & & & \gamma& \\ \hline
& -\alpha& & & & & & \\
1& & & & & & & \\
& & & -\alpha& & & & \\
& & -1& & & & & 
\end{array}
\right]\otimes Q_\psi, 
\end{equation*}
\begin{equation*}
\setlength{\arraycolsep}{1pt}
\renewcommand{\arraystretch}{1.1}
v_6 =
\left[
\begin{array}{cccc|cccc}
& & & & & & -\gamma\beta& \\
& & & & & & & -\gamma\beta \\
& & & & \gamma& & & \\
& & & & & -\gamma& & \\ \hline
& & -\beta& & & & & \\
& & & \beta& & & & \\
1& & & & & & & \\
& 1& & & & & & 
\end{array}
\right]\otimes Q_\psi, \quad
v_7 =
\left[
\begin{array}{cccc|cccc}
& & & & & & & \gamma\alpha\beta \\
& & & & & & \gamma\beta& \\
& & & & & -\gamma\alpha& & \\
& & & & \gamma& & & \\ \hline
& & & \alpha\beta& & & & \\
& & -\beta& & & & & \\
& \alpha& & & & & & \\
1& & & & & & & 
\end{array}
\right]\otimes Q_\psi.
\end{equation*}
By Lemma~\ref{lem:central-elements-regular} and the definition of $f_\psi$, for $c\in\{1_G, a, b, ab\}$, $z\in Z$, and $x\in\{0, 1\}$, we have $\ell_{(cz, x)}\star f_\psi = \psi(z)\bigl(\ell_{(c, x)}\star f_\psi\bigr)$.
Since the matrices $\ell_{(cz, x)}$ form a basis of $\CA$, the elements $f_\psi, v_1, v_2, \ldots, v_7$ span $\CA \star f_\psi$.
The product rules below follow from these matrix representations and the identity
\[
(A\otimes Q_\psi)\star(B\otimes Q_\psi) = (A\star B)\otimes Q_\psi \quad (A, B\in M_8(\C)), 
\]
which holds because $Q_\psi^2=Q_\psi$.
The element $f_\psi$ is the identity element of this ideal, and $v_i\star v_j = O \quad (i\neq j)$.
Their product is given in TABLE~\ref{tab:product-Afpsi}.
All blank entries in the table denote the zero matrix $O$.
\begin{table}[htbp]
\caption{The product table of $\CA \star f_\psi$}
\label{tab:product-Afpsi}
\centering
\renewcommand{\arraystretch}{1.5}
\begin{tabular}{c||c|ccccccc}
$\star$& $f_\psi$& $v_1$& $v_2$& $v_3$& $v_4$& $v_5$& $v_6$& $v_7$\\
\hline\hline
$f_\psi$& $f_\psi$& $v_1$& $v_2$& $v_3$& $v_4$& $v_5$& $v_6$& $v_7$\\
\hline
$v_1$ & $v_1$ & $\alpha f_\psi$ & & & & & & \\
$v_2$ & $v_2$ & & $\beta f_\psi$ & & & & & \\
$v_3$ & $v_3$ & & & $-\alpha\beta f_\psi$ & & & & \\
$v_4$ & $v_4$ & & & & $\gamma f_\psi$ & & & \\
$v_5$ & $v_5$ & & & & & $-\gamma\alpha f_\psi$ & &\\
$v_6$ & $v_6$ & & & & & & $-\gamma\beta f_\psi$ & \\
$v_7$ & $v_7$ & & & & & & & $\gamma\alpha\beta f_\psi$
\end{tabular}\end{table}

The product table also shows that these eight elements are linearly independent.
Indeed, suppose that
\[
c_0f_\psi+\sum_{i = 1}^{7}c_iv_i = O.
\]
For each $j$, write $v_j^{\star 2} = d_jf_\psi$, where $d_j\neq0$, and choose $k\neq j$.
Multiplying the relation by $v_j$ and then by $v_k$ gives $c_jd_jv_k = O$.
Since $v_k^{\star 2} = d_kf_\psi\neq O$, we have $v_k\neq O$, and hence $c_j = 0$.
Thus $c_j = 0$ for every $j$, and then $c_0 = 0$.
Consequently, $f_\psi, v_1, v_2, \ldots, v_7$ form a basis of $\CA \star f_\psi$.

Since all seven coefficients are nonzero, they define a nondegenerate symmetric bilinear form on a seven-dimensional vector space.
Therefore, 
\[
\CA \star f_{\psi} \cong \SPIN_{\C} (7).
\]
\textit{The case $\F = \R$.}
Suppose first that $\phi = \overline{\phi}$.
The same argument as in the complex case gives
\[
\CA \star e_{[\phi]} = \R e_\phi\cong\R.
\]

Suppose next that $\phi\neq\overline{\phi}$.
Regard $e_\phi$ and $e_{\overline{\phi}}$ as elements of $\CA_{\C}$, and set $u_{[\phi]} := \RI \left(e_\phi-e_{\overline{\phi}}\right)$.
Then $u_{[\phi]}\in\CA \star e_{[\phi]}$ and $u_{[\phi]}^{\star 2} = -e_{[\phi]}$.
Thus, $\{e_{[\phi]}, u_{[\phi]}\}$ is an $\R$-basis of $\CA \star e_{[\phi]}$, and the correspondence $e_{[\phi]}\longmapsto1, u_{[\phi]}\longmapsto \RI $ induces an isomorphism
\[
\CA \star e_{[\phi]}\cong\C.
\]
Let $\psi = \overline{\psi}$ and $\psi(s)\neq1$.
Then $\alpha, \beta, \gamma\in\{1, -1\}$.
For $1\leq i\leq7$, set $v_i := \ell_{q_i}\star f_{[\psi]}$.
As in the complex case, the elements $f_{[\psi]}, v_1, \ldots, v_7$ form a basis of $\CA \star f_{[\psi]}$, with $v_i\star v_j = O \quad (i\neq j)$.
The product rules are obtained from TABLE~\ref{tab:product-Afpsi}.
If $\alpha = \beta = \gamma = -1$, all seven diagonal coefficients are negative, and hence
\[
\CA \star f_{[\psi]} \cong \SPIN_{\R} (0, 7).
\]
In every other case, four diagonal coefficients are positive and three are negative.
Therefore, 
\[
\CA \star f_{[\psi]} \cong \SPIN_{\R} (4, 3).
\]
Finally, suppose that $\psi\neq\overline{\psi}$ and $\psi(s)\neq1$.
Regard $f_\psi$ and $f_{\overline{\psi}}$ as elements of $\CA_{\C}$, and set $u_{[\psi]} := \RI \left(f_\psi-f_{\overline{\psi}}\right)$.
Then $u_{[\psi]}^{\star 2} = -f_{[\psi]}$.
For $1\leq i\leq7$, define $v_i := \ell_{q_i}\star f_\psi$ and $v_i' := \ell_{q_i}\star f_{[\psi]}$.
The elements $f_{[\psi]}, v_1', \ldots, v_7', u_{[\psi]}, u_{[\psi]}\star v_1', \ldots, u_{[\psi]}\star v_7'$ form an $\R$-basis of $\CA \star f_{[\psi]}$.
The correspondence
\[
f_\psi\longmapsto f_{[\psi]}, \quad \RI f_\psi\longmapsto u_{[\psi]}, \quad v_i\longmapsto v_i', \quad \RI v_i\longmapsto u_{[\psi]}\star v_i'\quad (1\leq i\leq7)
\]
induces an isomorphism from the realification of $\CA_{\C} \star f_\psi$ onto $\CA \star f_{[\psi]}$.
Since $\CA_{\C} \star f_{\psi} \cong \SPIN_{\C} (7)$, we conclude that
\[
\CA \star f_{[\psi]} \cong \left( \SPIN_{\C} (7) \right)_{\R}.
\]
\end{proof}

\bibliographystyle{amsplain}
\providecommand{\bysame}{\leavevmode\hbox to3em{\hrulefill}\thinspace}
\providecommand{\MR}{\relax\ifhmode\unskip\space\fi MR }
\providecommand{\MRhref}[2]{%
  \href{http://www.ams.org/mathscinet-getitem?mr=#1}{#2}
}
\providecommand{\href}[2]{#2}

\providecommand{\bysame}{\leavevmode\hbox to3em{\hrulefill}\thinspace}

\end{document}